		\documentclass[12pt]{amsart}
\usepackage{amsmath,amsfonts,amssymb,amscd,amsthm,epsf}
	\newtheorem{prop}{Proposition}[section]
	\newtheorem{thm}[prop]{Theorem}
	\newtheorem{freedmanthm}[prop]{Freedman's Theorem}
	\newtheorem{Ethm}[prop]{Eliashberg's Theorem}
	\newtheorem{Mainthm}[prop]{Main Theorem}

	\newtheorem{questions}[prop]{Questions}

\theoremstyle{definition}

	\newtheorem{example}[prop]{Example}
	\newtheorem{examples}[prop]{Examples}

\theoremstyle{remark}
\def\ep{\varepsilon}
\def\X{{\chi}}
\def\complex{{\mathbb C}}
\def\PP{{\mathbb P}}
\def\CP{{\complex\PP}}
\def\oCP{\overline{\CP}}
\def\real{{\mathbb R}}
\def\zed{{\mathbb Z}}
\def\dim{\operatorname{dim}}
\def\Gr{\operatorname{Gr}}
\def\id{\operatorname{id}}
\def\Int{\operatorname{int}}
\def\mod{\operatorname{mod}}
\def\rel{\operatorname{rel}}
\begin{document}

\title{Stein surfaces as open subsets of $\complex^2$}
\author{Robert E. Gompf} 
\address{Department of Mathematics, The University of Texas at Austin,
1 University Station C1200, Austin, TX 78712-0257}
\email{gompf@math.utexas.edu}
\date{\today}
\thanks{Partially supported by NSF grant DMS-0102922.}
\begin{abstract}
An open subset $U$ of a complex surface can be topologically perturbed to 
yield an open subset whose inherited complex structure is Stein, if and 
only if $U$ is homeomorphic to the interior of a handlebody whose 
handles all have index $\le 2$.
\end{abstract}

\maketitle
\baselineskip=18pt

\section{Introduction}

Stein surfaces have been intensively studied  in some form 
by complex analysts since the early 20th century.
The most basic question, how plentiful Stein surfaces are, still remains 
wide open in various ways. 
For example, any open subset $U$ of a complex surface $X$ inherits a 
complex structure, so one may ask how commonly open subsets of $X$ satisfy 
the Stein condition.
This is already a deep classical problem when $X=\complex^2$. 
In the present paper, we show that Stein open subsets of a complex surface 
$X$ are ubiquitous in the following sense:
If $U$ satisfies the most basic topological condition necessary for the 
existence of a Stein structure, then after a suitable adjustment it becomes 
Stein as an open subset of $X$.
The adjustment is quite mild from a point-set topological viewpoint --- in 
particular, the homeomorphism type of $U$ is preserved, as well as the 
essential data of the topological embedding in $X$. 
However, the differential topology is radically altered. 
The diffeomorphism type of $U$ is typically changed, and frequently there 
are uncountably many possibilities for the resulting diffeomorphism type. 
The purpose of this paper is to state and prove the simplest form of the 
above existence result (Theorem~2.4), 
while developing the necessary background material. 
We also describe stronger results which will be presented in a forthcoming 
paper. 

\section{Stein surfaces and where to find them}

There are many equivalent ways of defining Stein surfaces. 
Perhaps the most efficient is that a {\em Stein surface} is a complex surface
(so the real dimension is 4) that admits a biholomorphic embedding as a 
closed subset of $\complex^N$ for some $N$. 
(A more relevant characterization will be given in Section~4.)
By the maximum modulus principle, Stein surfaces are never compact, but 
they might be considered the most natural generalization of compact 
K\"ahler surfaces. 
Global analysis such as Seiberg-Witten theory extends from the compact 
setting to the Stein setting \cite{KM}, because the Stein condition 
imposes ``pseudoconvexity'' at infinity, which has the effect of allowing 
suitable boundary conditions for the analysis. 
There is an analogous notion of convexity in symplectic geometry, leading 
to the concept of a Weinstein symplectic manifold \cite{W}, \cite{EG}, but 
for our present purposes the complex analytic setting is more natural. 
(We return to the Weinstein case in the last paragraph of the paper.)
It is important to note that, unlike affine algebraic surfaces, Stein 
surfaces need not have finite topology. 
Equivalently, a Stein surface need not be diffeomorphic to the interior of a 
compact 4-manifold with boundary. 
In fact, infinite topology is central to the construction in this paper, so
the resulting smooth manifolds are typically ``exotic.'' 
To apply Seiberg-Witten theory to such manifolds, for example, one must 
observe that they can be written as infinite nested unions of Stein surfaces 
with finite topology. 
(This at least gives information about compactly supported objects, 
for example the minimum genus of a representative of a fixed homology class.)

The most basic question we can can ask about Stein surfaces is how 
plentiful they are. 
We can make this question less vague in two different ways --- considering 
Stein surfaces as abstract manifolds or as open subsets of complex surfaces.

\begin{questions}\label{Q}
{\rm (a)} Which open, oriented 4-manifolds admit Stein structures?

{\rm (b)} In a fixed complex surface $X$, how common are open subsets $U$ 
that are Stein surfaces in the complex structure inherited from $X$?
\end{questions}

\noindent
The second question is already interesting for the simplest complex 
surface $X=\complex^2$. 
Such open subsets are classically called ``domains of holomorphy.'' 

To address these questions, it is convenient to use the language of 
handlebody theory.
Recall that a {\em $k$-handle} attached to the boundary of an $n$-manifold 
is a ``thickened up'' $k$-cell, that is, a copy of $D^k\times D^{n-k}$ 
smoothly attached along $\partial D^k \times D^{n-k}$, with the corners 
smoothed. 
The central $k$-cell $D^k \times\{0\}$ is called the {\em core} of the 
handle. 
A {\em handle decomposition} is a decomposition of a manifold into handles, 
essentially a thickened up $CW$-complex, cf.\ Figure~1. 
In the context of open manifolds, we build the handlebody first, then 
remove its boundary. 
Handlebody theory is equivalent to Morse theory:
A proper Morse function to $[0,\infty)$, such as height in Figure~1, 
determines a handle decomposition 
on its domain, with each index-$k$ critical point determining a $k$-handle 
(e.g.\ \cite{M}). 
This shows that all smooth manifolds admit handle decompositions. 

\begin{figure}
\centerline{\epsfbox{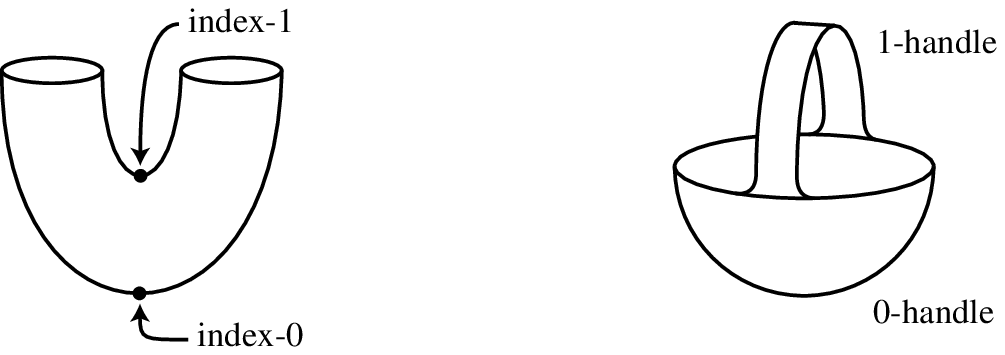}}
\label{handledecom}			
\caption{Handle decomposition of $S^1\times \real$}
\end{figure}

Question~\ref{Q}(a) was answered by Eliashberg \cite{E}, at least in principle.
His theorem is cleanest for higher dimensional {\em Stein manifolds} 
(complex $n$-manifolds embedding as closed subsets of $\complex^N$). 
It was classically known that every Stein manifold $X$ admits a proper Morse 
function to $[0,\infty)$ 
with critical points whose indices never exceed the complex  
dimension, or equivalently, a handle decomposition whose indices never 
exceed the middle (real) dimension. 
In addition, the complex structure on $X$ determines an almost-complex 
structure (complex vector bundle structure on $TX$). 
Surprisingly, these conditions are also sufficient for the existence of a 
Stein structure in high dimensions.

\begin{Ethm}\label{Yasha}
{\rm \cite{E}.} 
For $n\ge 3$, let $X$ be a smooth $2n$-manifold with an almost-complex 
structure $J$. 
If $X$ is the interior of a (possibly infinite) handlebody whose indices 
are all $\le n$, then $X$ admits a Stein structure homotopic to $J$.
\end{Ethm}

\noindent 
Since the existence of an almost-complex structure is a homotopy-theoretic 
question, Eliashberg's theorem characterizes manifolds admitting Stein 
structures in purely topological terms. 
The method also applies when $n=2$ 
(and the statement is well-known when $n=1$). 
In the $n=2$ case, almost-complex structures always exist (since we are 
dealing with oriented manifolds homotopy equivalent to 2-complexes). 
However, a framing obstruction arises, that always vanishes for $n\ge3$, 
but is nontrivial for 2-handles on 4-manifolds. 
Thus, Eliashberg characterizes oriented 4-manifolds admitting Stein structures 
as being interiors of handlebodies with all indices $\le2$ and an additional 
condition on the normal twisting of the 2-handle attaching maps. 
(See e.g.\ \cite{G1}  for a careful statement.) 
While this  is a purely topological characterization (and can be expressed 
entirely in terms of Kirby diagrams), it is somewhat 
difficult to apply in practice (e.g., \cite{G1}, \cite{GS}).

It is natural to ask whether the $n=2$ case of Theorem~\ref{Yasha} might 
be true as stated, even though it doesn't follow from Eliashberg's method. 
In fact, it already fails for the simplest example with a 2-handle: 
$S^2\times\real^2$ clearly satisfies the hypotheses of the theorem, being 
a 0-handle with a 2-handle attached in the simplest way. 
It is even complex $(\CP^1\times \complex)$. 
However, Seiberg-Witten theory \cite{LM}, \cite{KM} implies that any 
smoothly embedded, homologically essential sphere in a Stein surface must 
have homological self-intersection number $\le -2$, so $S^2\times\{0\} 
\subset S^2\times\real^2$ rules out the existence of a Stein structure. 
On the other hand, we might take inspiration from Freedman's revolutionary 
work \cite{F}, \cite{FQ}: 
While the main theorems of high-dimensional differential topology fail 
miserably in dimension~4, they tend to work if we ignore smooth structures 
and analyze the underlying topological manifolds. 
By combining Freedman's machinery with Eliashberg's, one obtains 
Theorem~\ref{Yasha} up to {\em homeomorphism}. 

\begin{thm}\label{G}
{\rm \cite{G1}.}
An oriented, topological 4-manifold $X$ is homeomorphic to a Stein surface 
if and only if it is the interior of a (possibly infinite) topological 
handlebody with all indices $\le 2$. 
In fact, every homotopy class of almost-complex structures on $X$ is 
then realized by such (orientation-preserving) homeomorphisms.
\end{thm}

\noindent
In dimension 4, topological handlebodies are uniquely smoothable (since 
we are gluing along 3-manifolds). 
See \cite{G1} for why almost-complex structures can be pulled back through 
homeomorphisms, up to homotopy. 
This theorem implies, for example, that there are Stein surfaces homeomorphic 
to $S^2\times \real^2$ (but of course not diffeomorphic to it), and in fact 
their Chern classes (hence, canonical classes) realize all even multiples 
of the generator of $H^2 (S^2 \times\real^2;\zed)$.

\begin{figure}[h]
\centerline{\epsfbox{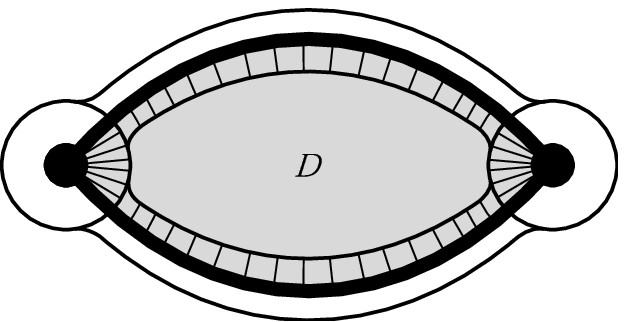}}
\label{fig2}
\caption{Core $D$ and its collar, comprising a 2-cell (shaded)}
\end{figure}

To address Question~\ref{Q}(b), we first observe that an open set $U\subset X$ 
has no chance of being Stein unless it has a handle decomposition with 
all indices $\le2$. 
Thus it suffices to consider open sets obtained by suitably thickening 
embedded 2-dimensional $CW$-complexes. 
For a topologically embedded 2-complex $K\subset X$, a neighborhood $U$ 
of $K$ will be called a {\em thickening} of $K$ if there is a homeomorphism 
$h:U\to \Int H$ onto the interior of a handlebody, sending cells to cores 
of handles. 
That is, there is a bijection between $k$-cells of $K$ and $k$-handles of $H$
for each $k$, and $h$ sends each cell to the core $D$ of the corresponding 
handle together with a collar of $\partial D$ that is radial with respect 
to the normal disks $\{p\}\times D^{4-\ell}$ of the lower-index handles. 
(See Figure~\ref{fig2}.)
A complex $K\subset X$ will be called {\em tame} if a thickening exists. 
This rules out wild embeddings analogous to  Alexander's horned sphere. 
Any two thickenings of a fixed $K\subset X$ will be homeomorphic 
($\rel K$ and preserving handles), although not necessarily diffeomorphic. 
We can ask whether they are topologically ambiently
isotopic (so they correspond under a homeomorphism $X\to X$ that is isotopic
to $\id_X$,
i.e., homotopic to $\id_X$ through homeomorphisms).
We may first need to topologically 
radially shrink the thickenings slightly to avoid technical difficulties at the 
boundary. 
(For example, the annuli $a<r<b$ in $\real^2$ are ambiently isotopic 
to each other for all $a,b>0$, but the case $a=0$ poses difficulties.) 
Then it seems plausible that the thickenings should be ambiently isotopic, 
but the question is apparently still open. 
Our main theorem now gives an answer to Question~\ref{Q}(b), by asserting 
that Stein open subsets $U\subset X$ are as common as possible, up to 
shrinking at the boundary and topological ambient isotopy. 
That is, the existence of any thickening implies (up to isotopy) the 
existence of a Stein thickening isotopic to it after shrinking.

\begin{Mainthm}\label{Main}
Let $K\subset X$ be a tame, topologically embedded $CW$ 2-complex in a 
complex surface, and choose $\ep>0$. 
Then after an $\ep$-small topological ambient isotopy of $K$, 
an $\ep$-neighborhood of $K$ contains a thickening $U$ of $K$ that is 
Stein in the complex structure that $U$ inherits as an open subset of $X$.
\end{Mainthm}

\noindent
By construction, the Stein thickening $U$ becomes topologically ambiently 
isotopic to the originally chosen thickening $V$, once we shrink $U$ and $V$ 
at their boundaries. 
A fancier construction \cite{G2} renders the shrinking unnecessary, provided 
that $\partial V$ is suitably tame. 
If $K$ is infinite, we may take $\ep$ to be a positive function on a 
neighborhood of $K$.

\begin{examples}\label{spheres}
(a) We consider $S^2$, the simplest complex of dimension~2. 
A smooth {\em 2-knot} is a smooth  embedding $S^2\hookrightarrow \real^4$. 
The theory of smooth 2-knots $K$ is even richer than classical knot theory. 
For example, all classical knot groups (as well as other groups) arise 
as smooth 2-knot groups $\pi_1(\real^4-K)$, and $\pi_2$ plays a role 
as well. 
For any 2-knot $K\subset \complex^2 = \real^4$, the theorem supplies a Stein
surface $U\subset \complex^2$ that remembers the original knot type. 
For example, $K$ and the closure of $U$ 
have complements with the same homotopy groups 
(although one must look carefully at the proof to see this, due to possible 
complications along the boundary of $U$). 
The Stein surface $U$ is homeomorphic to $S^2\times\real^2$ (since $K$ 
clearly has vanishing homological self-intersection number, so its normal 
bundle is trivial). 
Thus we obtain a plethora of counterexamples to the conjecture of \cite{Fo1} 
that no Stein open subset of $\complex^2$ is homotopy equivalent to $S^2$.
However, $U$  
is clearly not diffeomorphic to $S^2\times \real^2$, since $H_2(U)$ 
cannot be generated by a smoothly embedded sphere. 
(This shows the necessity of the topological isotopy --- $U$ cannot contain 
the original smooth sphere $K$ or anything smoothly isotopic to it, although 
by construction it contains a topological sphere topologically isotopic 
to $K$). 
One can arrange for the minimal genus of a smooth surface 
generating $H_2(U)$ to be any preassigned positive integer. 
Thus, each 2-knot $K$ corresponds to infinitely many diffeomorphism types 
of Stein surfaces $U$, distinguished by the minimal genus of the generator. 
On the other hand, one can identify a single diffeomorphism type that can be
realized by Stein thickenings (after topological isotopy)
of all smooth 2-knot types, and such universal 
diffeomorphism types can be found with arbitrary nonzero minimal genus.  
It is not known if topological 2-knots $S^2\hookrightarrow \complex^2$ 
exist that have topological $S^2\times\real^2$ thickenings but cannot be 
smoothed by a topological isotopy. 
If so, one would expect that no universal diffeomorphism types exist for 
their Stein thickenings, but any finite collection of topological 2-knots 
must (after isotopy) have a common diffeomorphism type of Stein thickenings.

(b) To get more complicated embeddings $S^2\hookrightarrow\complex^2$, 
start with a classical knot $\kappa$ (smooth embedding 
$S^1\hookrightarrow S^3$), and let $H_\kappa$ 
denote the handlebody obtained by attaching a 2-handle to a 0-handle $D^4$ 
along $\kappa\subset \partial D^4$, with the normal directions to $\kappa$ 
matched up so that $H_\kappa$ has trivial intersection pairing. 
The knot $\kappa$ is called (topologically or smoothly) {\em slice} if 
$H_\kappa$ embeds (topologically or smoothly) in $\real^4$ 
(or equivalently, in $S^4$). 
For example, if the Alexander polynomial $\Delta_\kappa(t)$ is 1, then 
$\kappa$ is topologically slice with $S^4-H_\kappa$ homotopy equivalent 
to $S^1$ \cite{FQ}. 
An embedding of $H_\kappa$ determines  an embedding $K$ of $S^2$ 
with a unique locally knotted point $p$, at which we see the cone on $\kappa$ 
in $D^4$. 
If we identify $K$ with the $CW$-complex $\{p\}\cup 2$-cell, then 
$\Int H_\kappa$ is a thickening of $K$. 
For a fixed topologically slice $\kappa$, we obtain an entire knot theory 
for such singular 2-spheres $K$. 
For example, if $\Delta_\kappa (t)=1$ we obtain topological knot complements in 
$\real^4$ realizing all the homotopy types of smooth 2-knot complements, 
by ambiently 
connected summing smooth 2-knots with the homotopy-$S^1$ case mentioned above.
Many knots with $\Delta_\kappa (t)=1$ are known not to be smoothly slice, 
so the corresponding singular topological 2-knots $K$ can never be 
isotoped to singular smooth 2-knots whose singularities are smooth cones. 
(That is, the embedded tame 2-complexes $K$ are unsmoothable.) 
However, the theorem still applies, and the discussion of (a) goes through 
without significant change, except that the resulting Stein surfaces 
in $\complex^2$ will all be homeomorphic to $\Int H_\kappa$ rather than to 
$D^2\times\real^2$, while they are still  homotopy equivalent to $S^2$.
(The minimal genus can still be taken to be arbitrarily large, 
although we lose precise control over it, and 
we still obtain universal diffeomorphism types in the smooth setting 
for a fixed, smoothly slice $\kappa$.) 
\end{examples}

The multitude of smooth structures we encounter here leads us to consider 
the quintessential phenomenon of open 4-manifold smoothing theory: 
exotic $\real^4$'s. 
These are smooth manifolds homeomorphic to $\real^4$ 
but not diffeomorphic to it. 
In contrast to dimensions $n\ne 4$, where exotic $\real^n$'s cannot exist, 
exotic $\real^4$'s realize uncountably many diffeomorphism types. 
(See Section~9.4 of \cite{GS} for a survey.) 
It was shown  in \cite{G1} that uncountably many exotic $\real^4$'s admit 
Stein structures. 
These smooth manifolds embed smoothly in the standard $\real^4$. 
Subsequently, L.~Taylor \cite{T} showed that uncountably many exotic 
$\real^4$'s (that do not embed in $\real^4$) require 3-handles in their 
handle decompositions, so they cannot support Stein structures. 
The techniques of the present paper can be used to prove the following 
theorem. 
(Details will appear in \cite{G2}.)

\begin{thm}\label{R4} 
There is a family of open subsets of $\complex^2$ (with compact closure) 
that are Stein and 
homeomorphic to $\real^4$, but realize uncountably many diffeomorphism 
types (with the cardinality of the continuum in ZFC set theory). 
\end{thm}

\noindent 
The usual way of finding an uncountable family of exotic $\real^4$'s is to 
identify one exotic $\real^4$ topologically with $\real^4$, then consider 
open balls of all sufficiently large radii. 
These topological balls inherit exotic smooth structures that can 
frequently be distinguished. 
The family in the above theorem arises in such a manner, with the Stein open 
subsets corresponding to radii lying in a Cantor set. 
Once we have such a family, we can connect it by a 1-handle to any $U$ 
as in Theorem~\ref{Main} and obtain a version of:

\begin{thm}\label{uncountable} 
For $K\subset \complex^2$ nonempty, 
the neighborhood $U$ given by Theorem~\ref{Main} 
lies in a family of Stein thickenings of $K$, 
nested (with compact closure in each other if $K$ is compact) 
with the order type of a Cantor 
set, and realizing uncountably many diffeomorphism types.
Every nonempty Stein open subset $U\subset \complex^2$ contains an 
uncountable nested family of Stein open subsets homeomorphic to $U$, 
but with no two diffeomorphic. 
\end{thm}

\noindent 
The structure of such families will be considered in more detail in 
Section~7 and in \cite{G2}. 
Similar (and sometimes stronger) results can be obtained for other cases 
of $K\subset X$. 
One would expect to essentially always obtain such uncountable collections  
of diffeomorphism types, but distinguishing these in complete generality 
seems beyond the range of present technology.

\begin{example}\label{spheres2}
For an embedded sphere $K\subset \complex^2$ as in Example~\ref{spheres}(a) 
or (b), we obtain an uncountable family of nested Stein surfaces corresponding 
to the given knot type, realizing uncountably many diffeomorphism types. 
This phenomenon is independent of the countably infinite families 
of the previous examples:
We can find a family as in Theorem~2.7 
for which the minimal genus of a generator is any fixed $g\ge1$ in (a) and 
fixed but arbitrarily large in (b), or arrange the minimal genus to increase 
without bound as the Stein surfaces get smaller. 
We can assume the intersection of the Stein 
surfaces is a sphere topologically ambiently isotopic to $K$, and that this 
is smooth except at a single point (aside from the smooth cone 
point $p$ in (b)).
\end{example}

Most of the remainder of this paper is devoted to a proof of 
Main Theorem~\ref{Main}, together with relevant background material. 
In the final Section~7, we will comment on the proofs of some of the 
remaining claims of this section, and present some additional observations.
More systematic and general statements and proofs will appear in \cite{G2}.

\section{Totally real surfaces}

Our first tool for proving Main Theorem~\ref{Main} involves embedded surfaces. 
Suppose $F$ is a compact, connected, oriented surface (so $\dim_\real F=2$)
smoothly embedded (or immersed) in a complex surface $X$. 
It is natural to ask how the complex structure of $X$ relates to $F$. 
We observe that there may be points $x\in F$ at which the 
(unoriented) tangent space 
$T_xF$ is a complex line; we call these {\em complex points}. 
Generically, $F$ has only finitely many complex points. 
This suggests counting complex points with suitable signs to obtain an index. 
We can define the sign of a complex point to be $(+)$ if and only if 
the complex orientation agrees with that of $F$.
However, there is also a second, more subtle, notion of sign, leading to the
distinction between elliptic and hyperbolic points (which should also be 
familiar to contact geometers), and resulting in two independent integer 
invariants (tracing back, e.g., to \cite{ChS}). 
While the details are not crucial for this paper, we state them for 
completeness: 
The Grassmann bundle $\Gr (TX)$ of oriented real 2-planes on $X$ has a pair 
$\PP_\pm (TX)$ of codimension-2 subbundles, consisting of complex lines, 
positively and negatively oriented.  
The bundle $TF$ can be interpreted as a section of $\Gr (TX)$ over $F$, so 
it has intersection numbers with $\PP_\pm(TX)$ yielding 
the required pair of integers. 

These invariants, repackaged by linear combination, can be recognized as 
familiar characteristic classes: 
the Chern number $\langle c_1(X),F\rangle$ and the sum $e(\nu F) +\chi (F)$ of 
normal and tangent Euler numbers. 
If $F$ is {\em totally real}, i.e., if $F$ has no complex points, then both 
invariants must vanish by definition. 
By a harder result of Eliashberg and Harlamov, these are the only obstructions.

\begin{thm}\label{real}
{\rm \cite{EH}.} 
A surface $F\subset X$ as above is smoothly isotopic to a totally real 
surface if and only if $\langle c_1 (X),F\rangle =0=e(\nu F) +\X (F)$.
\end{thm}

\noindent 
In fact, when both invariants vanish, Eliashberg and Harlamov connect 
canceling pairs of complex points by arcs in $F$, and then eliminate 
the complex points by a $C^0$-small twist of $F$ fixing each arc 
but rotating its normal directions. 
(See also \cite{N} for more details and \cite{Fo1} for a more 
general version.) 
The theorem is intended to apply rel boundary when $\partial F\ne \emptyset$. 
In that case, one needs suitable framings on $\partial F$ to interpret the 
characteristic classes. 
Generically, $F$ has no complex points on $\partial F$, so the tangent and 
outward normal vector fields $t$ and $n$ along $\partial F$ are a complex 
framing of $TX|\partial F$ for defining $\langle c_1 (TX),F\rangle$, and 
multiplying by $i$ sends $t$ to a normal vector field to $F$ along 
$\partial F$ defining the relative Euler class $e(\nu F)$. 
(Of course, the topological Euler characteristic $\X(F)$ is also the tangent 
Euler number relative to $t$ in this case.)

\section{Smooth Stein thickenings}

Our method for constructing Stein surfaces is taken from Eliashberg's 
paper \cite{E} proving Theorem~\ref{Yasha}. 
We will illustrate the method by attempting to prove a smooth version 
of Main Theorem~\ref{Main} --- so the open cells of $K$ are smoothly embedded, 
the isotopy is smooth, and the thickenings come from smoothly embedded 
handlebodies. 
We are doomed to fail, even for the simple case $K=S^2\subset \complex^2$,
but the method works for some pairs $K\subset X$. 
It is instructive to see how the method fails in general, so we can 
ultimately fix the proof by passing to the topological setting and invoking 
Freedman theory.

The Morse functions $\varphi$ discussed prior to 
Eliashberg's Theorem~\ref{Yasha} 
(with all indices $\le \dim_\complex X$) are actually {\em plurisubharmonic}
functions. 
These are essentially characterized as having level sets $\varphi^{-1}(c)$ 
that are all {\em pseudoconvex}. 
In complex dimension~2, this is the same as saying that the complex line 
field determined on each $\varphi^{-1}(c)$ (by the complex structure) is 
a contact structure that is positive in the boundary orientation on 
$\varphi^{-1}[0,c]$. 
(See, e.g., \cite{G1} or Chapter~11 of \cite{GS} for further discussion.) 
By a theorem of Grauert \cite{Gr}, a complex manifold with a proper 
plurisubharmonic function to $[0,\infty)$ must be Stein, 
so Stein manifolds are 
characterized by the existence of these functions. 
Eliashberg builds his Stein manifolds as complex manifolds with 
plurisubharmonic Morse functions, extending the complex structure and 
function one handle at a time. 
In the case of Stein open subsets (\cite{E}, Theorem~1.3.6), the complex 
structure is inherited from a carefully chosen embedding of the handle 
into the ambient manifold.

Eliashberg's method is simplest for 0- and 1-handles. 
In our case, we start with a smooth 2-complex $K\subset X$ in a complex 
surface. 
First thicken each 0-cell to a standard closed $\ep$-ball in some 
local holomorphic coordinates. 
Such a ball $B$ automatically has a pseudoconvex boundary that is a level 
set of the plurisubharmonic function $\varphi (z)=\|z\|^2$   
in the local chart near $B$. 
We can assume these balls are the 0-handles of our given smooth thickening 
of $K$, after a smooth isotopy of the thickening. 
Now each 1-cell of $K$ determines an embedding $\gamma :I= [0,1]\to X$ 
whose image intersects the 0-handles precisely at $\gamma (\partial I)$. 
After a $C^1$-small smooth isotopy fixing $\gamma$ and $d\gamma$ on 
$\partial I$, we can assume $\gamma$ is real-analytic.
Then complexifying $\gamma$ gives a holomorphic embedding from a neighborhood 
of $I$ in $\complex$. 
Using  a real-analytic vector field $v$ on $\gamma (I)$ normal to this complex 
curve, we can extend to a neighborhood of $I\times \{0\}$ in $\real^2$, 
then complexify to a neighborhood in $\complex^2$. 
This neighborhood of $I\times \{0\}$ in $\complex^2$ 
is Eliashberg's model complex 1-handle, on which 
he explicitly constructs a plurisubharmonic extension of $\varphi$. 
(Here we also need the 0-handle boundaries to intersect $I\times \{0\}$ 
orthogonally in the local model. 
This is easily achieved by suitable initial choice of $d\gamma$ and $v$ 
at $\partial I$.) 
After trimming each 1-handle back to a level set of $\varphi$, we have a 
Stein thickening $U_1$ of the 1-skeleton of $K$, and we can assume this 
agrees with the 0- and 1-handles of our original thickening (after an 
isotopy of the latter). 
This thickening $U_1$ has a natural boundary $\partial U_1$, which is a level
set of the plurisubharmonic function $\varphi$ defined near $U_1$. 
(If $K$ has infinitely many 0-cells, properness requires $\varphi$ to take 
arbitrarily large values on some 0-cells, so in practice we should add 
some 0-handles after we have begun adding 1-handles if we wish $\varphi$ 
to be proper. 
This causes no difficulties. 
We also add a thin collar at each stage to arrange $\varphi$ to be 
well-defined on the infinite union.  
Note that $K$ is at least locally finite since it embeds in a manifold.) 

Handles of larger index involve additional complications; for us, it 
suffices to consider a 2-handle. 
Each 2-cell of $K$ determines a smooth embedding $f:D^2\to X$ 
with $f^{-1} (\partial U_1) = \partial D^2$, which we would like 
to identify with the core of a model complex 2-handle, $D^2\times \{0\} 
\subset D^2\times D^2 \subset \real^2 \times i\real^2 \cong \complex^2$. 
This is only possible if $f(D^2)$ suitably resembles $D^2\times\{0\}$ 
in the model --- for example, it should be totally real.
(The model is chosen this way since the natural descending disks of a 
plurisubharmonic function are totally real. 
Further motivation is that every totally real submanifold $R$ has Stein 
tubular neighborhoods, since 
the squared distance to $R$ in any Hermitian metric is 
plurisubharmonic near $R$.
Note that in the previous paragraph, 
the 0- and 1-cells were automatically totally real.) 
First we adjust $\partial f(D^2)$. 
After a $C^1$-small perturbation, we can assume $\partial U_1$ is 
real-analytic near $\partial f(D^2)$. 
A $C^0$-small isotopy of $\partial f(D^2)$ in $\partial U_1$ makes 
it {\em Legendrian}, i.e., tangent to the contact planes of $\partial U_1$, 
and (\cite{E}, Lemma~2.5.1) a further $C^1$-small perturbation makes 
$f|\partial D^2$ also real-analytic. 
The Legendrian condition implies $f(D^2)$ has no complex points on its 
boundary, so the invariants of Theorem~\ref{real} are defined. 
If these obstructions vanish, then we can isotope $f$ $\rel\partial D^2$ 
to a totally real embedding. 
A $C^1$-small isotopy fixing $f$ and (if real-analytic) $df$ on 
$\partial D^2$ now makes $f$ real-analytic as well. 
(This follows from a general theorem, but for $D^2$ it can also be proved 
using Fourier series.) 
Now $f$ complexifies to a holomorphic embedding of a model 2-handle, 
and Eliashberg's construction gives the required plurisubharmonic extension 
of $\varphi$ as before. 
(The contact condition and suitable initial choice of outward normal 
to $f(D^2)$ along $\partial f(D^2)$ guarantee that $\partial U_1$ is 
orthogonal to $D^2\times \{0\}$ in the local model as required.)  

In conclusion, we can smoothly isotope $K$ and give it a Stein thickening, 
provided that the two integer obstructions of Theorem~\ref{real} vanish 
for each 2-cell of $K$. 
When can we arrange this? 
Note that the obstructions are  not well-defined until we fix a choice of 
the Legendrian curve $C= \partial f(D^2)$, since complex points can 
initially slide across the boundary. 
If we make one choice of $C$ and extend its tangent vector field $t$ to a 
section $v$ of the contact plane field near $C$
in $\partial U_1$, we can change to a different choice $C'$ 
by adding small spirals to $C$ parallel to the contact planes 
(winding relative to 
$v$), changing the frame $(t,n)$ used to define the relative 
characteristic classes. 
We can force $\langle c_1(X),F\rangle$ to vanish by this trick. 
(In fact, $\langle c_1(X), f(D^2)\rangle$ can be interpreted as 
the rotation number of the Legendrian curve $C$ relative to some $v$; 
it is well-known that this can be changed 
arbitrarily by adding zig-zags to the diagram.) 
The other invariant $e(\nu F) + \X (F) = e(\nu f(D^2))+1$ is more problematic. 
(Its vanishing translates to the condition that the normal framing 
must equal $tb(C)-1$, 
the extra requirement arising in the $n=2$ case of 
Eliashberg's Theorem~\ref{Yasha};  cf.\ \cite{G1} or \cite{GS}.) 
This invariant vanishes $\mod 2$ when the Chern obstruction does, since 
both count complex points with suitable signs. 
We can increase it by 2 without changing $\langle c_1(X),f(D^2)\rangle$, by 
adding an up/down pair of left-handed 
spirals (zig-zags), but cannot always decrease it.
This is the fundamental flaw of the method, which prevents us from solving 
cases like $S^2\subset \complex^2$. 
It is possible to solve some explicit cases, by directly computing the 
invariants and arranging them to vanish, but one cannot expect completely 
general results from this approach. 

\section{Casson handles}

We now introduce the necessary elements of topological 4-manifold theory, 
with historical motivation. 
By the 1960's, differential topology in high dimensions $(\ge 5)$ had become 
a mature field, thanks to powerful theorems for reducing topological 
questions to algebraic ones --- notably, the Surgery and $s$-Cobordism 
Theorems. 
The proofs of these theorems failed in dimension~4, and little was known 
about 4-manifolds before 1980. 
The failure of the proofs can be traced to a key step called the 
{\em Whitney trick}, which allows us to separate intersecting submanifolds 
of complementary dimension. 
Under suitable hypotheses, we can group the extra intersection points in 
pairs with opposite sign, then connect each pair by an arc in each submanifold, 
to obtain a circle $C$ as in Figure~3(a) . 

\begin{figure}[h]
\centerline{\epsfbox{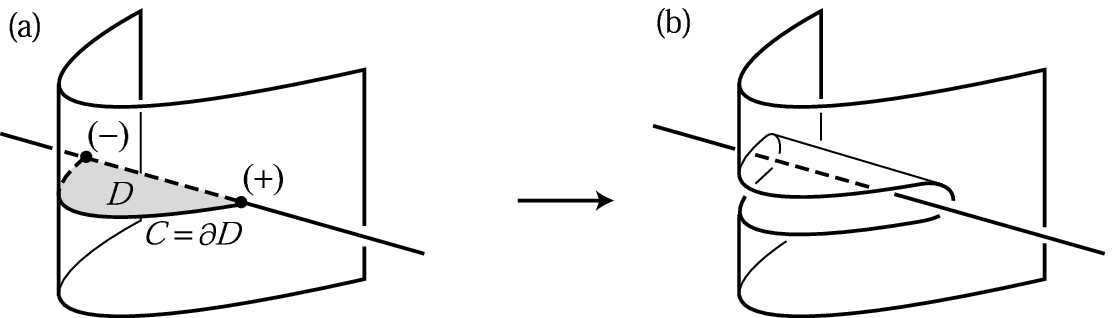}}
\label{fig3}
\caption{The Whitney trick} 
\end{figure}

\noindent
If we can find an embedded 2-disk $D$ with $\partial D=C$ and $\Int D$ 
disjoint from the submanifolds, then we can eliminate the pair of 
intersections by pushing one submanifold across $D$ as in Figure~3(b). 
(One must also pay attention to the directions normal to $D$, but that 
obstruction always vanishes in high dimensions.) 
In high dimensions, one can find the disk $D$ when needed: 
The hypotheses guarantee that $C$ is nullhomotopic in the complement of 
the submanifolds, and the resulting map of a disk can be assumed to be 
a smooth embedding by transversality $(2+2<5)$. 
In dimension~4, however, the best one can  obtain by transversality is a 
smoothly immersed disk with transverse double-point intersections. 
Thus the fundamental question of 4-manifold theory is when immersed disks 
can be replaced by embedded disks, or equivalently, embedded 2-handles. 

A major assault on this problem was launched by Casson in the 1970's 
\cite{C}. 
As he observed, a {\em kinky handle} $\kappa$, or closed tubular neighborhood 
of a generically immersed 2-disk in a 4-manifold, is 
most obviously distinguished from a 2-handle by its fundamental group: 
Each {\em kink} (self-intersection) contributes a generator to this free group.
One can recover a contractible space by attaching a 2-cell for each kink. 
In fact, Casson located a framed link $L$ in $\partial\kappa$, disjoint 
from the attaching region $\partial_-\kappa$, such that attaching 2-handles 
along $L$ transformed $(\kappa,\partial_-\kappa)$ into a 2-handle  
$(D^2\times D^2,S^1\times D^2)$. 
We can now rephrase the fundamental question: 
Given an embedding $(\kappa,\partial_-\kappa) \to (X,\partial X)$ (where we can 
take the 4-manifold $X$ to be the complement of a neighborhood of the 
surfaces we wish to separate by the Whitney trick), when can we find embedded
2-handles in $X - \Int \kappa$ attached along $L$? 
Such 2-handles would transform $\kappa$ into the desired 2-handle. 
Unfortunately, the new embedding problem seems even harder than the original.
While this may have stopped lesser mathematicians, Casson proceeded to show 
that, after suitable modification of the original immersion, one can at least 
find disjoint kinky handles attached to $L$. 
The resulting union with $\kappa$ is now called a 
{\em 2-stage Casson tower} $T_2$.
(One must use care with the meaning of attaching a kinky handle with a given 
framing, but this need not concern us here.) 
There is a new framed link on the boundary of the new kinky handles, on which 
attaching 2-handles would collapse the entire tower down to a 2-handle; 
Casson at least obtained disjoint kinky handles there for a 
{\em 3-stage Casson tower} $T_3$ in $(X,\partial X)$. 
By iterating the procedure, he obtained $n$-stage towers for all $n$, 
$\kappa = T_1 \subset T_2 \subset T_3 \subset \cdots \subset X$. 
The union of these towers, with the boundary removed except for 
$\Int \partial_-\kappa$, is now called a 
{\em Casson handle} $(CH,\partial_-CH)$.

Why should we care about such a complicated infinite construct?
First note that while each $T_n$ has a free fundamental group (generated 
by the top-stage kinks), a Casson handle is simply connected: 
Any loop $\gamma$ in $CH$ lies in some $T_n$ by compactness, but inclusion 
$T_n\hookrightarrow T_{n+1}$ is trivial in $\pi_1$. 
(The generators at each stage are killed by disks at the next stage.) 
Casson's more careful analysis showed that the smooth manifold 
$(CH,\partial_-CH)$ has the same proper homotopy type as an open 2-handle
$(D^2\times \real^2,S^1\times \real^2)$. 
At this point, Freedman picked up the project. 
After most of a decade studying the internal structure of Casson handles, 
and the application of some difficult point-set topology, Freedman 
proved his amazing theorem.

\begin{freedmanthm}\label{Freedman}
{\rm \cite{F}.}
Every Casson handle is homeomorphic to an open 2-handle.
\end{freedmanthm}

\noindent
Suddenly the Whitney trick and high-dimensional topology worked for 
topological 4-manifolds (first for trivial $\pi_1$, then later for many 
other fundamental groups \cite{FQ}). 
Immediate consequences included the 4-dimensional Poincar\'e Conjecture, 
a complete classification of closed, simply connected topological 
4-manifolds, and a Fields Medal for Freedman. 
(Some credit is also due for Quinn \cite{Q}, for strengthening the 
conclusion of the Classification Theorem.)
The main question left open, whether Casson handles might actually be 
diffeomorphic to open 2-handles, was answered negatively by Donaldson's 
contemporaneous work (yielding another Fields Medal): 
One cannot usually even find a smoothly embedded disk in a Casson handle, 
bounded by the attaching circle in $\partial_-CH$; otherwise the 
high-dimensional machinery would run smoothly, contradicting many 
results from gauge theory.
Fortunately, Donaldson's work did not precede Freedman's; otherwise 
topological 4-manifold theory might never have been developed. 

The Casson handle embedding theorem that we will need is taken from  
Quinn's paper (\cite{Q}, Proposition~2.2.4). 
A smooth 4-manifold $V$ with a circle $C$ in its boundary will be called 
a {\em possibly exotic 2-handle} with {\em attaching circle} $C$ if there 
is a homeomorphism $h:(D^2\times\real^2,S^1\times\{0\}) \to (V,C)$, and 
$h(D^2\times \{0\})\subset V$ will be called a {\em topological core}. 
(This core is unique up to topological ambient isotopy.) 
Since Casson handles are examples, we cannot in general expect to have a 
smoothly embedded 2-disk in $V$ bounded by $C$, and so $h$ cannot usually be 
smoothed near $D^2\times\{0\}$. 
However, we do obtain the following:

\begin{thm}\label{Quinn}
{\rm \cite{Q}.} 
A possibly exotic 2-handle $(V,C)$ always contains an unknotted Casson 
handle $CH$ attached along $C$.
\end{thm}

\noindent 
Here, {\em unknotted} means that the topological core of $CH$ given by 
Freedman's Theorem is topologically ambiently isotopic in $V$ 
$(\rel C)$ to the topological core $h(D^2\times \{0\})$ of $V$. 
Quinn had to settle for a ``weakly unknotted'' Casson handle; we prove 
the stronger statement using more modern technology. 
This theorem was originally assumed to be of limited use, since it only 
gives a Casson handle where we already have a topological 2-handle. 
Basically, it only replaces an arbitrary exotic 2-handle by a slightly 
better understood exotic 2-handle. 
However, we will find this upgrade crucial for creating Stein thickenings, 
and Quinn's proof will be an important part of our construction for 
Theorem~\ref{Main}.

\begin{proof} 
While the given homeomorphism  $h:D^2\times \real^2\to V$ cannot be 
smoothed near $D^2\times \{0\}$, Quinn  uses Freedman theory to 
investigate where $h$ {\em can} be smoothed. 
His Handle Straightening Theorem (2.2.2 of \cite{Q}) is based on 
showing that a 
homeomorphism can be smoothed  (isotoped to a diffeomorphism) 
outside of a small neighborhood of a smooth 2-complex with special properties. 
Thus a homeomorphism of an open 0- or 1-handle that is a diffeomorphism 
on the boundary can be smoothed near its core, 
simply by smoothly isotoping the core away from the 2-complex via 
transversality, whereas this fails for 2-handles. 
In the latter case, however, his theorem shows that we can avoid the 
problematic 2-complex by smoothly homotoping $D^2\times \{0\}$ $(\rel\partial)$ 
to a suitable generically immersed disk $D\subset D^2\times\real^2$. 
Then we can assume $h$ is a diffeomorphism near $D$, so $h(D)$ is a smoothly 
immersed disk in $V$ with $\partial D=C$. 

The advantage of using $h(D)$, rather than any generic smoothing of the map 
$h|D^2\times \{0\}$, is that the kinky handle $\kappa$ in $D^2\times\real^2$ 
obtained by thickening $D$ can be turned into a 2-stage tower with no 
kinks at the second stage --- and furthermore, this tower in $D^2\times\real^2$ 
is smoothly ambiently isotopic to a closed 
tubular neighborhood of $D^2\times\{0\}$ 
(so we preserve unknottedness). 
To see this, note that Quinn's homotopy consists of {\em finger moves} to 
push $D^2\times \{0\}$ off of the bad set. 
These moves consist of ``pushing the disk with one finger,'' along some arc 
$\gamma$ and back through itself as in Figure~4, 
creating an extra pair of intersections.
\begin{figure}
\centerline{\epsfbox{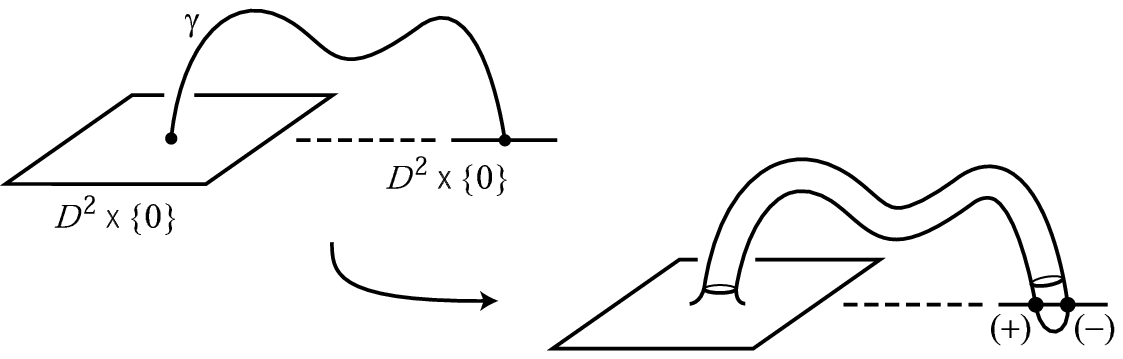}}
\label{fig4}
\caption{Finger move}
\end{figure}
Finger moves, which can be thought of as inverses of Whitney moves, have 
been used extensively for cleaning up immersed surfaces since Casson's work. 
A finger move is determined by the arc $\gamma$ with endpoints on the 
surface, together with a normal vector field to $\gamma$, whose 
orthogonal complement in $\nu\gamma$ carries the circle bundle comprising 
the tube around our finger. 
(Automorphisms of the circle bundle do not affect the image set $D$.) 
Since any two nonzero sections of the 3-plane bundle $\nu\gamma$ are homotopic 
rel boundary if they agree on the boundary, $D$ is determined by the 
arcs $\gamma_i$ specifying the finger moves. 
In fact, this collection of arcs is (in our application) smoothly 
isotopic to a standard model. 
This is because 1-manifolds cannot knot or link in a 4-manifold (homotopy 
implies isotopy) and the corresponding homotopy problem is trivial: 
Once we isotope $\gamma_i$ so that  it agrees near its endpoints 
with the corresponding model 
arc $\gamma'_i$, the two arcs determine an element of 
$\pi_1 (D^2\times\real^2 - D^2\times\{0\}) \cong\zed$. 
However, we can change this element by a generator, by turning $\gamma_i$ 
once around the normal fiber to $D^2\times\{0\}$ at one endpoint of $\gamma_i$.
Thus we can assume $\gamma_i$ and $\gamma'_i$ are homotopic, hence isotopic. 
Now that we have reduced to considering a standard model of a finger move, 
it is straightforward to exhibit the required pair of disks. 

We can now complete Quinn's construction. 
Applying the homeomorphism $h$ to our 2-stage tower in $D^2\times\real^2$, 
we obtain towers $h(\kappa) = T_1 \subset T_2^* \subset V$. 
While $T_1$ is a smooth 1-stage tower (kinky handle), our new disks in 
$D^2\times\real^2$ hit the non-smooth part of $h$, so the second stage 
of $T_2^*$ consists of topologically embedded 2-handles $h_i$. 
Stripping off the boundary of each $h_i$ (except for 
$\Int \partial_-h_i$, which is smooth), we obtain possibly exotic 
2-handles (in the smooth structure 
inherited from $V$), so we can repeat the previous construction in each $h_i$. 
Inductively assume we have towers $T_n\subset T_{n+1}^* \subset V$, where 
the $n$-stage tower $T_n$ is smooth, and $T_{n+1}^*$ is obtained from it 
by adding topologically embedded 2-handles. 
Then $T_{n+1}^*$ is homeomorphic to a 2-handle (since there are no top-stage
kinks); also assume $T_{n+1}^*$ is topologically ambiently isotopic to a 
closed tubular neighborhood of the topological core disk of $V$. 
Applying the previous construction to the possibly exotic 2-handles of the 
top stage, we recover the induction hypotheses with $T_n\subset T_{n+1} 
\subset T_{n+2}^* \subset T_{n+1}^*$. 
Taking the infinite union of $T_1\subset T_2\subset T_3\subset\cdots$ 
and removing the boundary except for $\Int \partial_-h (\kappa)$, we obtain 
a smoothly embedded Casson handle $CH\subset V$ with attaching circle $C$.

To prove $CH\subset V$ is unknotted, let $\Delta$ be a topological core disk 
of $CH$. 
Working in the topological category, we add boundary to $V$ to obtain 
$\Delta \subset D^2 \times D^2 = B^4$, with $V$ identified 
with $D^2\times \Int D^2$. 
We know that $(\Delta,\partial\Delta) \subset (B^4,\partial B^4)$ is 
{\em flat}, i.e., it is the core of an embedded 2-handle 
$(D^2\times D^2,S^1\times D^2) \hookrightarrow (B^4,\partial B^4)$ 
(obtained from the product structure of $CH\approx D^2\times\real^2$), 
whose interior we denote by $W\approx D^2\times\Int D^2$. 
We wish to show that $\pi_1 (B^4 -\Delta)\cong\zed$, for then Freedman's 
work immediately shows that $\Delta$ is unknotted in $B^4$ and hence in $V$:
In fact, $B^4-W$ is a topological 4-manifold with boundary $S^2\times S^1$ 
(since $\partial \Delta = C\subset \partial B^4$ is unknotted), 
so the $s$-Cobordism Theorem with $\pi_1\cong \zed$ 
(e.g., \cite{FQ}) identifies 
$B^4-W$ with $B^3\times S^1$, and this uniquely glues to $W$. 
Since $H_1(B^4-\Delta) \cong\zed$ (e.g., by Mayer-Vietoris), it suffices 
to show that any nullhomologous loop $\gamma$ in $B^4-\Delta$ is 
nullhomotopic. 
Compactness of $\Delta\subset CH$ guarantees that $\Delta$ lies in some 
subtower $T_n$ of $CH$. 
By construction, $T_n$ lies in $T_{n+1}^*$, which is a deformation retract 
of $B^4$ (since it is isotopic to a closed tubular neighborhood of the 
core of $V$). 
Thus, we can assume $\gamma$ lies in $T_{n+1}^* -\Delta$. 
Now we push $\gamma$ off of the top-stage 2-handles of $T_{n+1}^*$, so that 
$\gamma$ lies in $\Int T_n-\Delta\subset CH-\Delta$. 
But inclusion $CH-\Delta \hookrightarrow B^4-\Delta$ induces an 
$H_1$-isomorphism, since both groups are $\zed$ generated by a meridian of 
$\Delta$. 
Since $\gamma$ is nullhomologous in $B^4-\Delta$, it is then 
nullhomologous in $CH-\Delta$, and hence nullhomotopic (since 
$CH-\Delta\approx D^2\times (\real^2 -\{0\})$   has abelian $\pi_1$). 
\end{proof}

\section{Proof of Main Theorem~\ref{Main}}

We prove Theorem~\ref{Main} by combining the methods of the two previous 
sections. 
Recall that we are given a topologically embedded 2-complex $K$ in a complex 
surface $X$, and $K$ is tame, i.e., it can be thickened to a handlebody. 
We wish to move $K$ by an $\ep$-small topological ambient isotopy, so that 
it has an $\ep$-small Stein thickening. 
To make the construction $\ep$-small, we simply subdivide $K$ as 
necessary and work within the resulting small handles of a thickening.
The 1-skeleton of $K$ is still easy to deal with: 
By Quinn's Handle Straightening Theorem (2.2.2 of \cite{Q}), we can assume 
the 0-handles are embedded by homeomorphisms that are smooth near the 
vertices of $K$, then make them  Stein as in Section~4 (by a 
topological ambient isotopy fixing $K$ outside the 0-handles). 
Similarly, we can assume (after an $\ep$-small topological isotopy) 
that the 1-handles are 
smooth  near the cores, then modify them as in Section~4. 
We have now topologically isotoped $K$ so that the 1-skeleton is smooth 
and has a small Stein thickening $U_1$, with $\partial U_1$ a level set of 
a plurisubharmonic function $\varphi$, and $U_1$ extends to a topological 
thickening of all of $K$. 
(If $K$ is infinite, we again work with only finitely many handles at a 
time, starting on 2-handles before finishing the 1-handles, to maintain 
properness of $\varphi$.) 

The main challenge is the 2-cells.
Each 2-cell determines a topological embedding $f:D^2\to X$ that is 
smooth near $f^{-1}(\partial U_1)=\partial D^2$, and the image is the 
core of a possibly exotic 2-handle. 
The proof of Theorem~\ref{Quinn} transforms $f$ into a smooth 
immersion $f'$ such that the induced kinky handle has second stage 
topologically embedded 2-handles for reconstructing (up to topological 
isotopy) a thickening of $f(D^2)$. 
After we clean up near $\partial f'(D^2)$ as in Section~4, the 
obstructions of Theorem~\ref{real} are defined for $f'$. 
As before, we can assume the Chern class obstruction vanishes, so the 
remaining obstruction $e(\nu)+1$ is even. 
Recall that the fundamental difficulty of the method of Section~4 was the 
failure of this Euler class obstruction to vanish in general. 
However, we are no longer constrained to work with embedded surfaces. 
For immersed surfaces, we can change $e(\nu)$ by $\pm 2$, simply by adding a
kink with sign $\mp 1$. 
(Note that for a closed, oriented, generically immersed surface $F$, its 
homological self-intersection number is 
$F\cdot F = e(\nu F) + 2\text{ Self }F$, where $\text{Self }F$ is the 
signed number of double points of $F$. 
This follows by pushing off a copy $F'$ of $F$ and counting intersections 
with $F$ as in Figure~5 to obtain $F'\cdot F=F\cdot F$.  
\begin{figure}
\centerline{\epsfbox{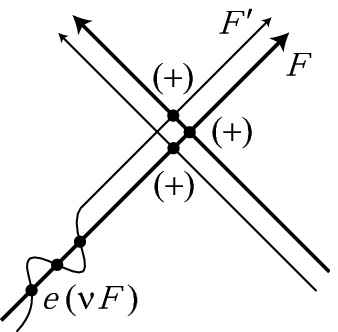}}
\label{fig5}
\caption{Homological self-intersection number}
\end{figure}
Since $F\cdot F$ is a homological invariant, adding a $\mp$ kink to $F$ 
must change $e(\nu F)$ by $\pm 2$.) 
Thus, we can make $e(\nu)+1$ vanish by adding kinks to $f'(D^2)$. 
Each kink can be added in a standard chart $(\real^4,\real^2)$ using a 
standard model that has a smoothly embedded second-stage 2-handle. 
Thus, the resulting immersion $f''$ still has the second stage 
required for Quinn's construction. 
We can easily adjust $f''$ near each double point so that the 
two sheets correspond to $\real^2$ and $i\real^2$ in a holomorphic local 
chart, then apply Theorem~\ref{real} away from the double points to make 
$f''$ totally real. 
A $C^1$-small perturbation, fixing $f''$ and $df''$ on $\partial D^2$ and 
at the double points, makes $f''$ real analytic and totally 
real, with multiplication by $i$ interchanging the tangent directions 
at each double point. 
Theorem~1.3.5 of \cite{E} says that totally real immersions, with 
$i$-invariant double points as above, have Stein regular neighborhoods, 
and its proof allows us to extend $\varphi$ over a kinky handle $\kappa$ 
with core $f''(D)$, so that $\partial \kappa-\partial_-\kappa$ 
lies in a level set.

To complete the construction, we apply induction as for Theorem~\ref{Quinn},
while at each stage adding a thin collar and 
maintaining the Stein condition as above. 
We proceed with the cells of $K$ and tower stages ordered by $\zed^+$ 
so that each stage of the construction is a finite handlebody, but each 
handle and kinky handle is eventually included. 
The union $U$ of these handlebodies is itself an open 
topological handlebody, whose 1-skeleton $\Int U_1$ is a thickening of the 
1-skeleton of $K$ (which we have already adjusted by a topological isotopy). 
The 2-handles of $U$ are Casson handles, one for each 2-cell of $K$. 
By the unknottedness assertion of Theorem~\ref{Quinn}, each  2-cell is 
topologically ambiently isotopic $\rel\partial$ to the core 
of the corresponding Casson handle, so after a topological ambient isotopy 
of $K$, $U$ becomes a topological thickening of $K$. 
(For infinite $K$, we can ensure that the ambient isotopy extends from the 
original thickening $V$  to all of $X$, by taking the 
preassigned $\ep$ to be a continuous
function with $\ep|V>0$ but $\ep |X-V=0$.)
Since $U$ was constructed with a proper plurisubharmonic function, it is 
Stein in the complex structure inherited from $X$, so it is the required 
thickening completing the proof of Theorem~\ref{Main}.\qed 

\begin{examples}\label{surfaces} 
a) For a smooth 2-knot $K$ in $\complex^2$, we have $c_1(\complex^2)=0$, 
$e(\nu K) = K\cdot K=0$ and $\X (K) =2$. 
Thus, we can make the sphere $K$ totally real by adding a single positive kink.
We can do this with a standard model kink so that the second stage disk
can also be made real with a single positive model kink.
By induction, one obtains a Stein thickening $U$ of $K$, after 
a $C^0$-small topological ambient isotopy that is smooth except near the 
kink, where $U$ is obtained by adding the simplest  Casson handle (with 
one positive kink at each stage) to a 4-ball along an unknot 
(with framing~0). 
This is a universal diffeomorphism type for Stein thickenings of smooth 
2-knots, in the sense of Example~\ref{spheres}, since we have just 
described $U$ up to diffeomorphism by data independent of $K$. 
(The diffeomorphism type of a Casson handle is determined by the number and
signs of the kinks of the kinky handles. 
This data can be expressed by a based tree with a sign attached to each 
edge.) 
The minimal genus of the generator of $H_2(U)$ is 1, since $U$ cannot 
contain a smooth essential sphere, but the immersed sphere can be smoothed 
to an embedded torus. 

b) More generally, for any smooth embedding $F\hookrightarrow X$ of a closed, 
oriented, connected surface into a (possibly noncompact) complex surface with 
$\langle c_1 (X),F\rangle =0$, we must add $m=\frac12\max (0,F\cdot F+ \X(F))$
positive kinks to the 2-cell to make it totally real. 
(Recall that instead of negative kinks, we can suitably adjust the 
Legendrian boundary when $F\cdot F+\X(F)<0$.) 
We obtain a Stein thickening where the Casson handle attached to the 
1-skeleton has only positive kinks, with $m$ at the first stage and one on 
each higher-stage kinky handle. 
(For $m=0$, this is just a smooth 2-handle, and the resulting Stein tubular 
neighborhoods can also be obtained by another method \cite{Fo1}, \cite{Fo2}.) 
This diffeomorphism type is then universal for all such embeddings of $F$
(with $F\cdot F$  fixed). 
When $\langle c_1 (X),F\rangle \ne 0$, we must absorb this 
obstruction with extra spirals 
in the Legendrian boundary, so the above holds with 
$m=\frac12\max (0,F\cdot F+\X(F) +|\langle c_1(X),F\rangle|)$. 
We obtain a universal diffeomorphism type $U_{g,k,c}$ of Stein 
thickenings for smooth embeddings with $g(F) =g$, $F\cdot F=k$ 
and $|\langle c_1(X),F\rangle| =c$. 
Note that we can realize any $g,c\ge0$ and $k \equiv c\mod 2$ by 
such embeddings, by immersing $F$ with degree~1 in a holomorphic line 
bundle realizing $g$ and $c$, then taking  $X$ to be a pulled back 
tubular neighborhood of $F$. 
The minimal genus of $U_{g,k,c}$ is $g +m$ (except conceivably when $g=0$, 
$m=1$ and $k\le -2$), 
by the genus bound of \cite{LM}. 
(See also \cite{GS} Theorem~11.4.7.) 
To realize universal diffeomorphism types of all larger minimal genera
(for fixed $g,k,c$), simply use $U_{g,k,c+2n}$ for all $n>0$. 
These can be substituted for $U_{g,k,c}$ by adding extra positive 
kinks at the first stage of the construction and compensating by suitable 
spirals in the attaching circle. 
For example, we obtain universal diffeomorphism types of any nonzero genus 
in Example~\ref{spheres}(a). 

c) Consider {\em topological} embeddings $F\hookrightarrow X$ that can be 
thickened (to a homeomorphism of the 2-plane bundle over $F$ with Euler 
number $F\cdot F$). 
Now our construction requires Quinn's Handle Straightening Theorem, so it is 
much  harder to extract information about the complexity of the resulting
Casson handle. 
For any fixed $g,k,c$ as above, we have topological embeddings 
$F\hookrightarrow X_n = U_{g,k,c+2n}$ for which the minimal genus 
realizing $[F]\in H_2(X_n)$ is arbitrarily large. 
Thus, there cannot 
be a universal diffeomorphism type for this problem. 
However, we will see that there is a fixed diffeomorphism type of $U$ for 
any finite collection of such embeddings with fixed $g,k$. 
\end{examples}

\section{More tricks}

The machinery of topological 4-manifold theory developed over a period of 
more than a decade, and many specialized techniques arose in the process.
One can obtain much more information about Stein thickenings by reviewing 
these techniques and inserting the word ``Stein'' in various places. 
We sketch some ideas here, and will return to the matter in \cite{G2}.

We first observe that any two Casson handles have a common refinement: 
Simply compare the Casson handles, from the first stage up, and add kinks 
of suitable sign to each whenever needed, to produce a new Casson handle 
smoothly embedded in both. 
(Equivalently, one constructs a signed tree containing both of the 
original signed trees.) 
Given a finite collection of topological embeddings 
$K\hookrightarrow X_i$, with 
homeomorphic $(\rel K)$ thickenings, we immediately see how to construct 
Stein thickenings for them that are diffeomorphic to each other:
For each 2-cell of $K$, construct the corresponding Casson handles in each 
$X_i$ simultaneously, adding extra kinks where necessary so that the 
resulting Casson handles all have the same signed tree, hence are 
diffeomorphic. 
(A sufficient excess of positive kinks everywhere guarantees that the 
constructed thickenings can be made Stein.) 
We can also enforce lower bounds on minimal genera by this method, refining 
the Stein thickening so that it smoothly embeds somewhere else where a genus 
bound for the homology class is known. 
For example, given $K\subset X$ as usual and $\alpha \ne 0\in H_2(K)$ 
(if this exists), Theorem~\ref{G} implies the thickening of $K$ is 
homeomorphic to a Stein surface $V$ for which $\langle c_1(V),\alpha\rangle$ 
exceeds any preassigned value, and hence the minimal genus of $\alpha$ in $V$
exceeds a preassigned value (\cite{GS} Exercise~11.4.11(d) and solution). 
A sufficiently refined Stein thickening $U$ for $K$ smoothly embeds in $V$
(not preserving $c_1$), so the minimal genus of $\alpha$ in $U$ also 
exceeds the preassigned value. 
This explains the genus bound of Example~\ref{spheres}(b). 

While smoothly embedded complexes provide a rich source of examples 
for Theorem~\ref{Main}, we would also like to locate nonsmooth examples. 
One good method is Casson's original theorem for embedding Casson 
handles \cite{C}. 
This shows, for example, that in a simply connected 4-manifold 
$X$, any class of square $\pm1$ in $H_2(X)$ whose orthogonal 
complement is odd can be represented by a topological embedding of 
$\pm \CP^2- \{p\}$ (hence, a tame topological $S^2$), and any hyperbolic pair 
$\left[\begin{smallmatrix} 0&1\\ 1&0\end{smallmatrix}\right]$ 
in the intersection form can be represented by an embedding of 
$S^2\times S^2 - \{p\}$ (so a tame $S^2\vee S^2$). 
(These can also be deduced from Freedman's Classification Theorem 
if $X$ is closed.) 
If $X$ is a closed, minimal complex surface with $b^+>1$,  gauge theory shows 
that these embeddings can never be made smooth, 
but Theorem~\ref{Main} makes the images Stein. 

Next we consider exotic $\real^4$'s. 
The basic constructions of these are in Casson's original paper, although 
he could not draw the final conclusions with neither Freedman nor 
Donaldson theory yet available.
The general class of exotic $\real^4$'s containing those of Theorem~\ref{R4} 
appeared in \cite{DF}, while the simplest case of these was shown to be 
exotic in \cite{BG} (see also \cite{GS}) and to be Stein in \cite{G1}. 
These simplest examples are constructed with one 0-handle, two 1-handles, 
one 2-handle and a Casson handle. 
The range of diffeomorphism types comes from varying the complexity of the 
Casson handle. 
To embed such an example as a Stein open subset of $\complex^2$, we must 
first embed the honest handles. 
While it is easy to embed the 1-skeleton as a Stein surface, it takes work 
to arrange the 2-cell so that the obstructions of Theorem~\ref{real} can 
both be explicitly seen to vanish. 
(Details will appear in \cite{G2}.) 
Once this is achieved, the Casson handle can be  located by Casson's 
Embedding Theorem, and our previous techniques yield the Stein surface.

One of the main ingredients of the proof of Freedman's Theorem~\ref{Freedman}
is his Reimbedding Theorems.
These elucidate the structure of a Casson handle $CH$ by finding other 
Casson handles with compact closure inside. 
Ultimately, Freedman obtains an uncountable family of Casson handles 
indexed by a standard Cantor set $C\subset I= [0,1]$, nested with compact 
closure inside each other with the ordering inherited from $C$. 
Instead of a Casson handle at $0\in C$, Freedman shows that the intersection 
of the remaining Casson handles 
can be taken to be a topological core disk $\Delta$ 
of $CH$, and $\Delta$ can be assumed smooth except at a single interior point.
The main principle of our Section~6 is that in a complex surface, Casson 
handles can be assumed to preserve Stein structures, as long as we are 
allowed to add enough positive kinks. 
We conclude that any Stein surface $U$ 
made with finitely many handles and Casson handles 
contains a family of Stein surfaces 
homeomorphic to $U$, nested with compact closure, 
with the order type of $C-\{0\}$. 
(For example, any $U$ arising from a finite handlebody in 
Theorem~\ref{G} or a compact $K$ in Theorem~\ref{Main} has this structure.)
The intersection of these nested Stein surfaces 
will be a tame 2-complex, smoothly embedded 
except at a single point in the interior of each 2-cell, and the 2-cells 
will be totally real where they are smooth. 
In the case of the above exotic $\real^4$, \cite{DF} shows that each 
diffeomorphism type appears only countably often within the nested family, 
so there are uncountably many diffeomorphism types of Stein exotic $\real^4$'s 
realized this way in $\complex^2$, with the cardinality 
of the continuum in ZFC set theory (proving Theorem~\ref{R4}).  
The argument of \cite{DF} is too specialized to apply to all Stein 
surfaces, but in some cases (e.g., any open subset of $\complex^2$) we 
can distinguish uncountably many diffeomorphism types after we connect 
by a 1-handle to the family of Theorem~\ref{R4} 
(proving Theorem~\ref{uncountable}). 
Other methods sometimes distinguish {\em all} diffeomorphism types in a nested 
family as above, e.g., for Stein surfaces homeomorphic to 
$\pm \CP^2 -\{p\}$ or $S^2\times S^2-\{p\}$. 
Our nesting procedure can be combined with genus restrictions, so that 
the minimal genus of a class $\alpha\ne 0\in H_2(U)$ increases without 
bound as the depth in the nesting increases.
Alternatively, we can arrange the minimal genus of $\alpha$ to be constant, 
by locating a surface $F$ realizing $\alpha$ with minimal genus in $U$, 
observing 
that for some $n$, $F$ lies in the union of the 1-skeleton and $n$-stage 
towers of $U$, then applying reimbedding only to the Casson handles attached 
on top of the $n$-stage towers.
Note that neighborhood systems of homotopy-equivalent Stein surfaces are 
of interest to complex analysts (e.g., \cite{Fo1}, \cite{Fo2} for 
neighborhood systems of smoothly embedded surfaces), but such neighborhoods 
are more typically diffeomorphic to each other. 

The family of exotic $\real^4$'s in \cite{DF} is actually parametrized 
by $I$, with $C\subset I$ corresponding to exotic $\real^4$'s made with 
reimbedding technology. 
To interpolate to a connected family in this way, one must know that the 
point-set boundaries of the nested $\real^4$'s are homeomorphic to $S^3$
(concentric round spheres in the topological $\real^4$ structure).
Unfortunately, embedded Casson handles typically have point-set 
boundaries that are not manifolds. 
(We can at least arrange that the complements of their closures 
have the same homotopy groups as the complements of their cores, by 
sharpening the method of proof of Theorem~\ref{Quinn}, 
but otherwise the boundaries are ugly.) 
To fix this, one replaces Casson handles by the infinite towers of capped 
gropes in \cite{FQ}. 
This means that between successive stages of immersed disks, one inserts 
many stages of embedded surfaces. 
This cleans up the boundaries, so that the topological $D^2\times \Int D^2$ 
can be parametrized with the radial coordinates corresponding to 
$C\subset I$ realized by reimbedded towers. 
We can apply this technology in the Stein setting, since increasing 
genus has the same effect on $e(\nu F) +\X(F)$ as adding positive kinks, 
but there are more technical details. 
(For example, the towers of \cite{FQ} are constructed, for convenience, 
with the same number of positive and negative kinks in each kinky handle, 
whereas we need an excess of positive kinks.) 
The final result is that the families of Stein surfaces 
parametrized by a Cantor set lie in continuous 1-parameter families 
of homeomorphic open sets (not necessarily Stein for parameter values in 
$I-C$). 

Our final observation on Stein surfaces is that Theorem~\ref{G} can be 
deduced from Main Theorem~\ref{Main}. 
(Of course, their proofs rely on similar ideas in any case.) 
If $X$ is the interior of a handlebody whose indices are all $\le2$, 
and $J$ is an almost-complex structure on $X$, there is a map 
$f$ from $X$ to the rational surface $\CP^2 \,\#\, \oCP^2$,  covered 
by a complex bundle map of the tangent bundles. 
(Send the 1-skeleton to a point, and wrap each 2-cell around the 
generators as necessary to preserve $c_1$.) 
By immersion theory, we can homotope $f$ to an immersion, then pull back 
the complex structure, to create a complex structure on $X$ homotopic to $J$.
If $K\subset X$ is the 2-complex obtained from the cores of the handles, 
applying Theorem~\ref{Main} to $K$ gives the required Stein surface.

Since Stein surfaces are in some sense equivalent to Weinstein 4-manifolds 
in the symplectic category, it is natural to ask whether there is an 
analog of Theorem~\ref{Main} for identifying open subsets of symplectic 
4-manifolds that are Weinstein. 
The main difference is that cores of 2-handles are Lagrangian in the 
Weinstein setting --- a closed condition that is strictly stronger than 
the totally real condition for Stein 2-handles. 
Fortunately, there is an $h$-principle for immersed Lagrangian 
submanifolds, so we can still obtain immersed disks by killing the invariants 
of Theorem~\ref{real}. 
It is not presently clear to the author, however, whether these immersions 
can be created with sufficient control to allow Quinn's tower construction 
to proceed. 
The author hopes to return to this question soon. 


\end{document}